\documentclass[12pt,leqno]{article}
\usepackage{amsopn,amsmath,amssymb,amscd,latexsym}
\usepackage[all]{xy}
\newcommand{\Z}{{\mathbb{Z}}}
\newcommand{\oI}{\overline{I}}
\newcommand{\id}{\mathrm{id}}
\newcommand{\spec}{\mathrm{Spec}\,}
\newcommand{\Aut}{\mathrm{Aut}\,}
\newcommand{\uAut}{\underline{\Aut}}
\newcommand{\coker}{\mathrm{coker}\,}
\newcommand{\End}{\mathrm{End}\,}
\newcommand{\uEnd}{\underline{\End}}
\newcommand{\Gl}{\mathrm{Gl}\,}
\newcommand{\Fh}{{\mathcal F}}
\newcommand{\Gh}{{\cal G}}
\newcommand{\Ih}{{\mathcal I}}
\newcommand{\Oh}{{\mathcal O}}
\newcommand{\emm}{{\mathfrak{m}}}
\newcommand{\oemm}{\overline{\emm}}
\newcommand{\tei}{\, | \,}
\newtheorem{theorem}{Theorem}
\newtheorem{lemma}{Lemma}
\begin{document}

%

\title{Representability of $\uAut_{\Fh}$ and $\uEnd_{\Fh}$}
\author{Niko Naumann}
\date{\today }
\maketitle

\begin{abstract}
  Recently N. Nitsure showed that for a coherent sheaf $\Fh$ on a noetherian scheme the automorphism functor $\uAut_{\Fh}$ is representable if and only if $\Fh$ is locally free. Here we remove the noetherian hypothesis and show that the same result holds for the endomorphism functor $\uEnd_{\Fh}$ even if one asks for representability by an algebraic space.
\end{abstract}


\section{Statement of results}\label{s1}
\subsection{}\label{s11}
  Let $X$ be a scheme and $\Fh$ a quasi-coherent $\Oh_X$-module of finite presentation. We are interested in the representability of the following two functors on the category of $X$-schemes:
\begin{eqnarray*}
\uAut_{\Fh} (X') & := & \Aut_{\Oh_{X'}} (f^* \Fh) \\
\uEnd_{\Fh} (X') & := & \End_{\Oh_{X'}} (f^* \Fh)
\end{eqnarray*}
for $f : X' \to X$ an $X$-scheme.

The result is as follows:

\begin{theorem}
\label{thm11}
Let $X$ be a scheme and $\Fh$ a quasi-coherent $\Oh_X$-module of finite presentation. Then the following are equivalent:
  \begin{enumerate}
  \item [1)] $\Fh$ is locally free.
  \item [2)] $\uAut_{\Fh}$ is representable by a scheme.
  \item [2')] $\uEnd_{\Fh}$ is representable by a scheme.
  \end{enumerate}
If $X$ is locally noetherian, these conditions are also equivalent to the following:
\begin{enumerate}
\item [3)] $\uAut_{\Fh}$ is representable by an algebraic space.
\item [3')] $\uEnd_{\Fh}$ is representable by an algebraic space.
\end{enumerate}
\end{theorem}

\subsection{}\label{s12}
 The equivalence of 1) and 2) in theorem \ref{thm11} in case $X$ is noetherian is the main result of \cite{N}. Our proof follows the ideas of {\em loc.cit.} closely. The main steps are contained in the following two lemmas:

\begin{lemma}
  \label{t1}
Let $A$ be a local ring and $M$ a finitely presented $A$-module which is {\em not} free. Then there is a local homomorphism $A \to B$ such that
\[
M \otimes_A B \cong B^n \oplus (B/b)^m \; ,
\]
for some $0 \neq b \in B , b^2 = 0$ and $m \ge 1,n \ge 0$.\\
If $A$ is noetherian, $B$ can be chosen to be artin.
\end{lemma}

We observe that in the last statement of the lemma the noetherian hypothesis is indispensable: let $(B , \emm)$ be a local ring such that there is $0 \neq b \in \bigcap_{n \ge 1} \emm^n$. Clearly $(b^2) \subsetneq (b)$, so after dividing out $(b^2)$ one gets a ring $B$ as in the lemma but for any local homomorphism $f : B \to C$ with $C$ {\em noetherian} one clearly has $f (b) = 0$.

\begin{lemma}
  \label{t2}
Let $S$ be a scheme and $S_0 \subseteq S$ a closed subscheme defined by a nilpotent ideal sheaf. Assume $X$ is a flat $S$-scheme and $f : X \to Y$ is an $S$-morphism such that $f \times \id_{S_0}$ is an isomorphism. Then $f$ is an isomorphism.
\end{lemma}

\subsection{} \label{s13}
 In order to treat the representability of $\uEnd_{\Fh}$ we will use the following observation:

\begin{lemma}
  \label{t3}
Under the assumptions of 1.1 the obvious natural transformation of (set-valued) functors $\uAut_{\Fh} \to \uEnd_{\Fh}$ is relatively representable by an open immersion.
\end{lemma}

 For completeness we also include a proof of the next lemma which is essentially lemma 5 of \cite{N} and shows the relative representability of a ``parabolic'' sub-group functor:\\
Let $X$ be a scheme and
\begin{equation}
  \label{eq:1}
  0 \longrightarrow \Fh' \longrightarrow \Fh \longrightarrow \Fh'' \longrightarrow 0
\end{equation}
a short exact sequence of quasi-coherent $\Oh_X$-modules with $\Fh'$ finitely presented and $\Fh''$ locally free. For any morphism $f : Y \to X$, the sequence $f^* ((\ref{eq:1}))$ is exact because $\Fh''$ is in particular $\Oh_X$-flat and it makes sense to consider
\[
P (Y) := \{ \alpha \in \Aut_{\Oh_Y} (f^* \Fh) \tei \alpha (f^* \Fh') \subseteq f^* \Fh' \} \subseteq \uAut_{\Fh} (Y) \; .
\]

\begin{lemma}
  \label{t4}
In the above situation, the natural transformation $P \hookrightarrow \uAut_{\Fh}$ is relatively representable by a closed immersion.
\end{lemma}

For basic facts about (relative) representability we refer to \cite{BLR}, 7.6.

\section{Proofs}
\subsection{}
\label{s21}
 In this subsection we dispense with the easy implications of theorem \ref{thm11}, the assumptions and notations of which we now assume:\\
As $\uAut_{\Fh}$ and $\uEnd_{\Fh}$ are clearly Zariski sheaves the problem of representing them is Zariski local on $X$, i.e. we can assume that $X$ is affine and $\Fh$ corresponds to a free module of finite rank. In this case, representability of both $\uAut_{\Fh}$ and $\uEnd_{\Fh}$ is obvious; we have proved the implications 1) $\Rightarrow$ 2) and 1) $\Rightarrow$ 2'). Finally, the implications 2) $\Rightarrow$ 3) and 2') $\Rightarrow$ 3') are trivial.

\subsection{}
\label{s22}
 
 {\bf Proof of lemma \ref{t1}:} Let $(A , \emm)$ be a local ring and $M$ a finitely presented $A$-module which is not free. We will find the required local homomorphism $A \to B$ as a suitable quotient of $A$:\\
Let
\begin{equation}
  \label{eq:2}
  A^m \xrightarrow{\alpha} A^n \xrightarrow{\beta} M \rightarrow 0
\end{equation}
be a minimal presentation of  $M$, i.e. $n = \dim_k (M / \emm M)$ where $k := A / \emm$ is the residue field of $A$. Then $M$ is free if and only if $\alpha = 0$: clearly $\alpha = 0$ is sufficient for freeness of $M$ and conversely, if $M$ is free, it is necessarily so of rank $n$, hence $\beta$ is a surjective endomorphism of $A^n$ which must be an isomorphism by a standard application of Nakayama's lemma, c.f. \cite{M}, thm. 2.4., hence $\alpha = 0$.\\
For any $J \subseteq \emm$, (\ref{eq:2}) $\otimes_A A / J$ is a minimal presentation of the $A / J$-module $M / JM$. If we denote by $I \subseteq A$ the ideal generated by the coefficients of any matrix representation of $\alpha$ and note that the minimality of (\ref{eq:2}) implies $I\subseteq \emm$ we find that $M / JM$ is $A / J$-free if and only if $\alpha\otimes id_{A/J}=0$ if and only if $I \subseteq J$. As $M$ is not $A$-free we have $I \neq 0$ and as $I$ is finitely generated we get $\emm I \subsetneq I$, again by Nakayama's lemma. By Zorn's lemma, using again that $I$ is finitely generated, there is an ideal $J$ with $\emm I \subseteq J \subsetneq I$ and which is maximal subject to these conditions (indeed, any ascending chain of such ideals admits its union as an upper bound because I is finitely generated). 

We claim that $B := A / J$ is as required:\\
By the maximality of $J$ the ideal $\oI := I / J$ is non-zero principal: $\oI = (b) , 0 \neq b \in B$ and we neccessarily have $b^2 = 0$: if not, we would have $b \in (b^2)$, i.e. $b = xb^2$ or $b (1-xb) = 0$ for some $x \in B$. As $b \in \oemm := \emm / J$, the maximal ideal of $B ,  1-xb$ was a unit of $B$, so we would have $b = 0$.

We now show that $M\otimes_{A} B$ has the desired structure: any coefficient $\alpha_{ij}$ of a matrix representation of $\alpha \otimes \id_B$ is of the form $\alpha_{ij} = bu_{ij} , u_{ij} \in B$. As by construction $\oemm b = 0$ we see that if $\alpha_{ij} \neq 0$, then $u_{ij} \in B^*$. We get a matrix equation $(\alpha_{ij}) = b (u_{ij})$ and $(u_{ij})$ can be chosen with $u_{ij} = 0$ or $u_{ij} \in B^*$, all $i,j$. Then the usual Gau{\ss}-algorithm can be applied to $(u_{ij})$, showing that indeed $M \otimes_A B \cong B^n \oplus (B / b)^m$ for some $m,n\geq 0$. As, by construction, $M \otimes_A B$ is not $B$-free, we finally see that $m \ge 1$. 

If $A$ is noetherian we can start the construction of $B$ by first dividing out a suitable high power of $\emm$: Indeed, if $M / \emm^n M$ was free for all $n \ge 1$ we would have $I \subseteq \bigcap_{n \ge 1} \emm^n = (0)$. Then the ring $B$ we obtain in the above construction is noetherian local with $\oemm$ nilpotent, hence zero-dimensional, i.e. $B$ is artin local.

{\bf Proof of lemma \ref{t2}:} We can assume that the ideal sheaf $\Ih$ of $S_0 \subseteq S$ satisfies $\Ih ^2 = 0$. Our assertion is local on $S, X$ and $Y$ and thus reduces to the following:\\
Given a ring $k$ and an ideal $I \subseteq k$ of square zero, if $f : A \to B$ is a morphism of $k$-algebras with $B$ $k$-flat and such that $f \otimes_k \id_{k / I}$ is an isomorphism, then $f$ is an isomorphism:\\
1) $f$ is surjective: any $b \in B$ can be written
\[
b = f (a) + \sum_j\alpha_jb_j'\; \mbox{ ,some} \; \alpha_j \in I , b'_j \in B , a \in A \; .
\]
Applying this to the $b'_j$ we get (for some $\alpha_{ij} \in I , b''_{ij} \in B , a_j \in A)$:
\[
b = f (a) + \sum_j\alpha_j(f(a_j)+\sum_{ij}\alpha_{ij}b_{ij}'')=f(a+\sum_j\alpha_j a_j )\; .
\]
2) $f$ is injective: For $K := $ker$ (f)$ the $k$-flatness of $B$ implies $K / IK = 0$ and the same argument as in 1) shows that $K = 0$.

  {\bf Proof of 2) $\Rightarrow$ 1) in theorem \ref{thm11}:} Under the notations of \ref{s11} we assume that $\uAut_{\Fh}$ is representable by a scheme and, by contradiction, that $\Fh$ is not locally free. Note that the assumption on representability is stable under base-change $Y\rightarrow X$. So, base-changing to a suitable local ring of $X$, we can assume $X = \spec (A)$ with $A$ local and $\Fh$ corresponding to a finitely presented $A$-module $M$ which is not free. According to lemma \ref{t1} we can assume $M \cong A^n \oplus (A / a)^m$ for some $0 \neq a \in A$ with $a^2 = 0$ and $m \ge 1$. Let $G \to S:= \spec (A)$ be the group-scheme representing $\uAut_M$ and put $S_0 := \spec (A / a)$. The sub-functor $G' \hookrightarrow G$ of automorphisms preserving (base-changes of) the direct summand $(A / a)^m$ is represented by a closed sub-group scheme (still to be denoted $G'$) according to lemma \ref{t4}.

Let $P \subseteq \Gl_{n+m , S}$ denote the standard parabolic sub-group of automorphisms preserving the rank $m$ direct summand. $P$ is flat over $S$, as can be seen over $\spec (\Z)$. There is a morphism of $S$-groups $f : P \to G'$ which on points is given by sending $\left( 
  \begin{smallmatrix}
    \alpha & 0 \\ \beta & \gamma
  \end{smallmatrix} \right)$ to $\left( 
  \begin{smallmatrix}
    \alpha & 0 \\ \pi \beta & \overline{\gamma}
  \end{smallmatrix} \right)$, where $\alpha \in \Aut_A (A^n) , \gamma \in \Aut_A (A^m) , \beta : A^n \to A^m$ and $\overline{\gamma} \in \Aut_{A/a} ((A / a)^m)$ denotes the reduction of $\gamma$ and $\pi : A^m \to (A / a)^m$ is the natural map. This ``point-wise'' description of $f$ is immediately checked to be functorial and a homomorphism and hence does indeed define a morphism of $S$-groups. Obviously, $f \times \id_{S_0}$ is an isomorphism, hence so is $f$ by lemma \ref{t2}. This is however a contradiction, because $f (S) : P (S) \to G' (S)$ is not injective, as $f (S) (\id_{A^n} \oplus (1-a) \id_{A^m}) = 1$ and $a \neq 0 , m \ge 1$.

\subsection{}
\label{s23} 
 {\bf Proof of lemma \ref{t3}:} Given a scheme $X$, a quasi-coherent $\Oh_X$-module $\Fh$ of finite presentation and some $\varphi \in \End_{\Oh_X} (\Fh)$ we have to show that there is an open sub-scheme $X_0 \subseteq X$ such for all $f : Y \to X , f^* (\varphi) \in \Aut_{\Oh_Y} (f^* \Fh) \subseteq \End_{\Oh_Y} (f^* \Fh)$ if and only if $f$ factors through $X_0$. Consider $\Gh := \coker (\varphi)$ and the exact sequence of $\Oh_X$-modules
    \begin{equation}
      \label{eq:3}
      \Fh \xrightarrow{\varphi} \Fh \xrightarrow{} \Gh \xrightarrow{} 0 \; .
    \end{equation}
We claim that $f^* (\varphi)$ is an automorphism if and only if $f^* (\Gh) = 0$: as $f^* ((\ref{eq:3}))$ is exact, necessity is obvious. If, conversely, $f^* (\Gh) = 0$ then for any $y \in Y$ $f^* (\varphi)_y$ is a surjective endomorphism of the finitely generated $\Oh_{Y,y}$-module $\Fh_y$, hence is an isomorphism, hence so is $f^* (\varphi)$.

So the sought for $X_0 \subseteq X$ is the complement of the support of $\Gh$ which is open, because $\Gh$ is finitely presented.

  {\bf Proof of lemma \ref{t4}:} Given a scheme $X$ and a short exact sequence $0 \to \Fh' \to \Fh \to \Fh'' \to 0$ of quasi-coherent $\Oh_X$-modules with $\Fh'$ finitely presented and $\Fh''$ locally free and some $\alpha \in \Aut_{\Oh_X} (\Fh)$, we have to show the representability by a closed sub-scheme of $X$ of the following functor on $X$-schemes:
\[
F (Y \xrightarrow{f} X) := \left\{ 
  \begin{array}{ccl}
* & , & f^* (\alpha) (f^* \Fh') \subseteq f^* \Fh' \\
\emptyset &, & \mbox{otherwise} \; .
  \end{array} \right.
\]
Clearly, $F$ is a Zariski sheaf, so the problem is local on $X$, i.e. we can assume that $X = \spec (A)$ is affine, $\Fh''$ corresponds to some $A^n$, $\Fh'$ corresponds to some $A$-module $M$ for which there is a presentation $A^a \to A^b \to M \to 0$ and $\Fh$ corresponds to some $A$-module $N$. The exact sequence $0 \to \Fh' \to \Fh \to \Fh'' \to 0$ then becomes an exact sequence $0 \to M \xrightarrow{\iota} N \xrightarrow{\pi} A^n \to 0$ of $A$-modules and we are given some $\alpha \in \Aut_A (N)$. Consider $\nu := \pi \alpha \iota$. As all the above sequences are exact after {\it any} base-change, we have $F (Y \xrightarrow{f} X) \neq \emptyset \iff f^* (\nu) = 0$.

We have a diagram (defining $\psi$):
\[
\xymatrix{
A^a \ar[r] & A^b \ar[r] \ar[dr]^{\psi} & M \ar[r] \ar[d]^{\nu} & 0 \\
           &                           & A^n                   &
}
\]
which is exact after any base-change, hence $f^* (\nu) = 0 \iff f^* (\psi) = 0$, for any $f : Y \to X$. So the closed sub-scheme of $X$ we are looking for is the one defined by the ideal of $A$ generated by the coefficients of any matrix representation of $\psi$.

  {\bf Proof of 2') $\Rightarrow$ 1) in theorem \ref{thm11}:} Under the notations of \ref{s11} we assume that $\uEnd_{\Fh}$ is representable by a scheme. Then so is $\uAut_{\Fh}$ by lemma \ref{t3}, hence $\Fh$ is locally free by what has been shown in \ref{s22}.

\subsection{}\label{s24}
   \quad \\
    {\bf Proof of 3) $\Rightarrow$ 1) and 3') $\Rightarrow$ 1) in theorem \ref{thm11}:} Under the notations of 1.1 we assume that $X$ is locally noetherian and that either 3) or 3') holds as well as, by contradiction, that $\Fh$ is not locally free. By lemma \ref{t3} we know in either case that $\uAut_{\Fh}$ is representable by an algebraic space. Using the last assertion of lemma \ref{t1} we can assume $X = \spec (A)$ with $A$ artin local. Then $\uAut_{\Fh}$ is representable by a scheme according to \cite{K}, p. 25, 7) contradicting what we proved in \ref{s22}.

{\bf Acknowledgements.} I would like to thank H. Frommer and M. Volkov for interesting discussions and G. Weckermann for excellent type-setting.

{\small Mathematisches Institut der WWU M\"unster\\
Einsteinstr. 62\\
48149 M\"unster\\
Germany\\
e-mail: naumannn@uni-muenster.de}

\end{document}